\pdfoutput=1
\RequirePackage{ifpdf}
\ifpdf % We are running pdfTeX in pdf mode
\documentclass[pdftex]{sigma}
\else
\documentclass{sigma}
\fi

\numberwithin{equation}{section}
\newtheorem{Theorem}{Theorem}[section]
\newtheorem{Corollary}[Theorem]{Corollary}
\newtheorem{Lemma}[Theorem]{Lemma}
\newtheorem{Proposition}[Theorem]{Proposition}
{\theoremstyle{definition}
\newtheorem{Example}[Theorem]{Example}
\newtheorem{Remark}[Theorem]{Remark}
}
\DeclareMathOperator{\Aut}{Aut}
\DeclareMathOperator{\charop}{char}
\DeclareMathOperator{\Der}{Der}
\DeclareMathOperator{\Hom}{Hom}
\DeclareMathOperator{\Gr}{Gr}
\DeclareMathOperator{\Id}{Id}
\DeclareMathOperator{\Lie}{Lie}
\DeclareMathOperator{\rmM}{M}
\DeclareMathOperator{\TwDer}{TwDer}

\begin{document}

\newcommand{\arXivNumber}{1312.4018}

\allowdisplaybreaks

\renewcommand{\PaperNumber}{065}

\FirstPageHeading

\ShortArticleName{Matched Pair Deformations and the Factorization Index}

\ArticleName{Bicrossed Products, Matched Pair Deformations\\ and the Factorization Index for Lie Algebras}

\Author{Ana-Loredana AGORE~$^{\dag\ddag}$ and Gigel MILITARU~$^\S$}

\AuthorNameForHeading{A.L.~Agore and G.~Militaru}

\Address{$^\dag$~Faculty of Engineering, Vrije Universiteit Brussel, Pleinlaan 2, B-1050 Brussels, Belgium}
\EmailD{\href{mailto:ana.agore@vub.ac.be}{ana.agore@vub.ac.be}, \href{mailto:ana.agore@gmail.com}{ana.agore@gmail.com}}
\URLaddressD{\url{http://homepages.vub.ac.be/~aagore/}}

\Address{$^\ddag$~Department of Applied Mathematics, Bucharest University of Economic Studies,\\
\hphantom{$^\ddag$}~Piata Romana 6, RO-010374 Bucharest 1, Romania}

\Address{$^\S$~Faculty of Mathematics and Computer Science, University of Bucharest,\\
\hphantom{$^\S$}~Str.
Academiei 14, RO-010014 Bucharest 1, Romania}
\EmailD{\href{mailto:gigel.militaru@fmi.unibuc.ro}{gigel.militaru@fmi.unibuc.ro},
\href{mailto:gigel.militaru@gmail.com}{gigel.militaru@gmail.com}}
\URLaddressD{\url{http://fmi.unibuc.ro/ro/departamente/matematica/militaru_gigel/}}

\ArticleDates{Received January 20, 2014, in f\/inal form June 10, 2014; Published online June 16, 2014}

\Abstract{For a~perfect Lie algebra $\mathfrak{h}$ we classify all Lie algebras containing $\mathfrak{h}$ as a~subalgebra of
codimension~$1$.
The automorphism groups of such Lie algebras are fully determined as subgroups of the semidirect product $\mathfrak{h} \ltimes
(k^* \times \Aut_{\Lie} (\mathfrak{h}))$.
In the non-perfect case the classif\/ication of these Lie algebras is a~dif\/f\/icult task.
Let $\mathfrak{l} (2n+1, k)$ be the Lie algebra with the bracket $[E_i, G] = E_i$, $[G, F_i] = F_i$, for all $i = 1, \dots, n$.
We explicitly describe all Lie algebras containing $\mathfrak{l} (2n+1, k)$ as a~subalgebra of codimension $1$ by computing all
possible bicrossed products $k \bowtie \mathfrak{l} (2n+1, k)$.
They are parameterized by a~set of matrices ${\rmM}_n (k)^4 \times k^{2n+2}$ which are explicitly determined.
Several matched pair deformations of $\mathfrak{l} (2n+1, k)$ are described in order to compute the factorization index of some
extensions of the type $k \subset k \bowtie \mathfrak{l} (2n+1, k)$.
We provide an example of such extension having an inf\/inite factorization index.}


\Classification{17B05; 17B55; 17B56}

\section{Introduction}

The theory of Lie algebras is among the most developed f\/ields in algebra due to his broad applicability in dif\/ferential
geometry, theoretical physics, quantum f\/ield theory, classical or quantum mechanics and others.
Besides the purely algebraic interest in this problem, the classif\/ication of Lie algebras of a~given dimension is a~central
theme of study in modern group analysis of dif\/ferential equations~-- for further explanations and an historical background
see~\cite{popovici}.
The Levi--Malcev theorem reduces the classif\/ication of all f\/inite-dimensional Lie algebras over a~f\/ield of characteristic zero
to the following three subsequent problems: (1)~the classif\/ication of all semi-simple Lie algebras (solved by Cartan); (2)~the
classif\/ication of all solvable Lie algebras (which is known up to dimension~$6$~\cite{gra2}) and (3)~the classif\/ication of all
Lie algebras that are direct sums of semi-simple Lie algebras and solvable Lie algebras.

Surprisingly, among these three problems, the last one is the least studied and the most dif\/f\/icult.
Only in 1990 Majid~\cite[Theo\-rem~4.1]{majid} and independently Lu and Weinstein~\cite[Theo\-rem~3.9]{LW} introduced the
concept of a~\emph{matched pair} between two Lie algebras $\mathfrak{g}$ and $\mathfrak{h}$.
To any matched pair of Lie algebras we can associate a~new Lie algebra $\mathfrak{g} \bowtie \mathfrak{h}$ called the
\emph{bicrossed product} (also called \emph{double Lie algebra} in~\cite[Definition~3.3]{LW}, \emph{double cross sum}
in~\cite[Proposition~8.3.2]{majid2} or \emph{knit product} in~\cite{Mic}).
In light of this new concept, problem~(3) can be equivalently restated as follows: for a~given (semi-simple) Lie algebra
$\mathfrak{g}$ and a~given (solvable) Lie algebra $\mathfrak{h}$, describe the set of all possible matched pairs $(\mathfrak{g},
\mathfrak{h}, \triangleleft, \triangleright)$ and classify up to an isomorphism all associated bicrossed products $\mathfrak{g}
\bowtie \mathfrak{h}$.
Leaving aside the semi-simple/solvable case this is just the \emph{factorization problem} for Lie algebras~-- we refer
to~\cite{abm1} for more details and additional references on the factorization problem at the level of groups, Hopf algebras,
etc.

The present paper continues our recent work~\cite{am-2013a, am-2013b} related to the above question (3), in its general form,
namely the factorization problem and its converse, called the \emph{classifying complement problem}, which consist of the
following question: let $\mathfrak{g} \subset \mathfrak{L}$ be a~given Lie subalgebra of $\mathfrak{L}$.
If a~complement of $\mathfrak{g}$ in $\mathfrak{L}$ exists (that is a~Lie subalgebra $\mathfrak{h}$ such that $\mathfrak{L} =
\mathfrak{g} + \mathfrak{h}$ and $\mathfrak{g} \cap \mathfrak{h} = \{0\}$), describe explicitly, classify all complements and
compute the cardinal of the isomorphism classes of all complements (which will be called the \emph{factorization index}
$[\mathfrak{L}: \mathfrak{g}]^f$ of $\mathfrak{g}$ in~$\mathfrak{L}$).
Our starting point is~\cite[Proposition~4.4]{am-2013b} which describes all Lie algebras~$\mathfrak{L}$ that contain a~given Lie
algebra~$\mathfrak{h}$ as a~subalgebra of codimension~$1$ over an arbitrary f\/ield~$k$: the set of all such Lie algebras~$L$ is
parameterized by the space~$\TwDer (\mathfrak{h})$ of twisted derivations of~$\mathfrak{h}$.
The pioneer work on this subject was performed by K.H.~Hofmann:~\cite[Theorem~I]{Hopmann} describes the structure
of~$n$-dimensional real Lie algebras containing a~given subalgebra of dimension~$n-1$.
Equivalently, this proves that the set of all matched pairs of Lie algebras $(k_0, \mathfrak{h}, \triangleleft, \triangleright)$
(by $k_0$ we will denote the Abelian Lie algebra of dimension~$1$) and the space $\TwDer (\mathfrak{h})$ of all twisted
derivations of $\mathfrak{h}$ are in one-to-one correspondence; moreover, any Lie algebra~$\mathfrak{L}$ containing
$\mathfrak{h}$ as a~subalgebra of codimension~$1$ is isomorphic to a~bicrossed product $k_0 \bowtie \mathfrak{h} =
\mathfrak{h}_{(\lambda, \Delta)}$, for some $(\lambda, \Delta) \in \TwDer (\mathfrak{h})$.
The classif\/ication up to an isomorphism of all bicrossed products $\mathfrak{h}_{(\lambda, \Delta)}$ is given in the case when~$\mathfrak{h}$ is perfect.
As an application of our approach, the group $\Aut_{\Lie} (\mathfrak{h}_{(\lambda, \Delta)}) $ of all automorphisms of
such Lie algebras is fully described in Corollary~\ref{izoaut}: it appears as a~subgroup of a~certain semidirect product
$\mathfrak{h} \ltimes(k^* \times \Aut_{\Lie} (\mathfrak{h}))$ of groups.
At this point we mention that the classif\/ication of automorphisms groups of all indecomposable real Lie algebras of dimension up
to f\/ive was obtained recently in~\cite{fisher} where the importance of this subject in mathematical physics is highlighted.
For the special case of sympathetic Lie algebras $\mathfrak{h}$, Corollary~\ref{morfismeperfectea} proves that, up to an
isomorphism, there exists only one Lie algebra that contains $\mathfrak{h}$ as a~Lie subalgebra of codimension one, namely the
direct product $k_0 \times \mathfrak{h}$ and $\Aut_{\Lie} (k_0 \times \mathfrak{h}) \cong k^* \times \Aut_{\Lie}
(\mathfrak{h})$.
Now, $k_0$ is a~subalgebra of $k_0 \bowtie \mathfrak{h} = \mathfrak{h}_{(\lambda, \Delta)}$ having $\mathfrak{h}$ as
a~complement: for a~$5$-dimensional perfect Lie algebra all complements of $k_0$ in $\mathfrak{h}_{(\lambda, \Delta)}$ are
described in Example~\ref{prefectnecomplet} as matched pair deformations of $\mathfrak{h}$.
Section~\ref{mpdefo} treats the same problem for a~given $(2n+1)$-dimensional non-perfect Lie algebra $\mathfrak{h}:=\mathfrak{l}(2n+1,k)$.
Theorem~\ref{teorema11} describes explicitly all Lie algebras containing $\mathfrak{l} (2n+1, k)$ as a~subalgebra of codimension~$1$.
They are parameterized by a~set ${\mathcal T}(n)$ of matrices $(A, B, C, D, \lambda_0, \delta) \in {\rmM}_n (k)^4 \times k
\times k^{2n+1}$: there are four such families of Lie algebras if the characteristic of~$k$ is $\neq 2$ and two families in
characteristic~$2$.
All complements of $k_0$ in two such bicrossed products $k_0 \bowtie \mathfrak{l} (2n+1, k)$ are described by computing all
matched pair deformations of the Lie algebra $\mathfrak{l} (2n+1, k)$ in Propositions~\ref{mpdef} and~\ref{mpdef2}.
In particular, in Example~\ref{n1} we construct an example where the factorization index of $k_0$ in the $4$-dimensional Lie
algebra $\mathfrak{m} (4, k)$ is inf\/inite: that is $k_0$ has an inf\/inite family of non-isomorphic complements in $\mathfrak{m}(4,k)$.
To conclude, there are three reasons for which we considered the Lie algebra $\mathfrak{l} (2n+1, k)$ in Section~\ref{mpdefo}:
on the one hand it provided us with an example of a~f\/inite-dimensional Lie algebra extension $\mathfrak{g} \subset \mathfrak{L}$
such that $\mathfrak{g}$ has inf\/initely many non-isomorphic complements as a~Lie subalgebra in $\mathfrak{L}$.
On the other hand, the Lie algebra $\mathfrak{l} (2n+1, k)$ serves for constructing two counterexamples in Remark~\ref{2noi}
which show that some properties of Lie algebras are not preserves by the matched pair deformation.
Finally, having~\cite[Corollary 3.2]{am-2012} as a~source of inspiration we believe that any $(2n+1)$-dimensional Lie algebra is
isomorphic to an~$r$-deformation of $\mathfrak{l} (2n+1, k)$ associated to a~given matched pair: a~more general open question is
stated at the end of the paper.

\newpage

\section{Preliminaries}
%\label{prel}
All vector spaces, Lie algebras, linear or bilinear maps are over an arbitrary f\/ield~$k$.
The Abelian Lie algebra of dimension~$n$ will be denoted by $k^n_0$.
For two given Lie algebras $\mathfrak{g}$ and $\mathfrak{h}$ we denote by $\Aut_{\Lie} (\mathfrak{g})$ the group of
automorphisms of $\mathfrak{g}$ and by $\Hom_{\Lie} (\mathfrak{g}, \mathfrak{h})$ the space of all Lie algebra maps
between $\mathfrak{g}$ and $\mathfrak{h}$.
A~Lie algebra $\mathfrak{L}$ \emph{factorizes} through $\mathfrak{g}$ and $\mathfrak{h}$ if $\mathfrak{g}$ and $\mathfrak{h}$
are Lie subalgebras of $\mathfrak{L}$ such that $\mathfrak{L} = \mathfrak{g} + \mathfrak{h}$ and $\mathfrak{g} \cap \mathfrak{h}
= \{0\}$.
In this case $\mathfrak{h}$ is called a~\emph{complement} of $\mathfrak{g}$ in $\mathfrak{L}$; if $\mathfrak{g}$ is an ideal of
$\mathfrak{L}$, then a~complement $\mathfrak{h}$, if it exists, is unique being isomorphic to the quotient Lie algebra
$\mathfrak{L}/ \mathfrak{g}$.
In general, if $\mathfrak{g}$ is only a~subalgebra of $\mathfrak{L}$, then we are very far from having unique complements; for
a~given extension $\mathfrak{g} \subset \mathfrak{L}$ of Lie algebras, the number of types of isomorphisms of all complements of
$\mathfrak{g}$ in $\mathfrak{L}$ is called the \emph{factorization index} of $\mathfrak{g}$ in $\mathfrak{L}$ and is denoted~by
$[\mathfrak{L}: \mathfrak{g}]^f$~-- a~theoretical formula for computing $[\mathfrak{L}: \mathfrak{g}]^f$ is given
in~\cite[Theorem 4.5]{am-2013a}.
For basic concepts and unexplained notions on Lie algebras we refer to~\cite{EW, H}.

A \emph{matched pair} of Lie algebras~\cite{LW, majid2} is a~system $(\mathfrak{g}, \mathfrak{h}, \triangleleft,
\triangleright)$ consisting of two Lie algebras~$\mathfrak{g}$ and~$\mathfrak{h}$ and two bilinear maps $\triangleright:
\mathfrak{h} \times \mathfrak{g} \to \mathfrak{g}$, $\triangleleft: \mathfrak{h} \times \mathfrak{g} \to \mathfrak{h}$ such that
$(\mathfrak{g}, \triangleright)$ is a~left $\mathfrak{h}$-module, $(\mathfrak{h}, \triangleleft)$ is a~right
$\mathfrak{g}$-module and the following compatibilities hold for all~$g, h \in \mathfrak{g}$ and~$x, y \in \mathfrak{h}$
\begin{gather*}
x \triangleright [g, h]  =  [x \triangleright g, h] + [g, x \triangleright h] + (x \triangleleft g) \triangleright h - (x
\triangleleft h) \triangleright g,
%\label{mpLie1}
\\
\left[x, y \right] \triangleleft g  =  \left[x, y \triangleleft g \right] + \left[x \triangleleft g, y \right] + x \triangleleft
(y \triangleright g) - y \triangleleft (x \triangleright g).
%\label{mpLie2}
\end{gather*}

Let $(\mathfrak{g}, \mathfrak{h}, \triangleleft, \triangleright)$ be a~matched pair of Lie algebras.
Then $\mathfrak{g} \bowtie \mathfrak{h}: = \mathfrak{g} \times \mathfrak{h}$, as a~vector space, is a~Lie algebra with the
bracket def\/ined~by
\begin{gather*}
%\label{bicrossed}
\{(g, x), (h, y) \}:= \big([g, h] + x\triangleright h - y \triangleright g, [x, y] + x \triangleleft h - y \triangleleft g\big)
\end{gather*}
for all~$g, h \in \mathfrak{g}$ and~$x, y \in \mathfrak{h}$, called the \emph{bicrossed product} associated to the matched
pair $(\mathfrak{g}, \mathfrak{h}, \triangleleft, \triangleright)$.
Any bicrossed product $\mathfrak{g} \bowtie \mathfrak{h}$ factorizes through $ \mathfrak{g} = \mathfrak{g} \times \{0\}$ and $
\mathfrak{h} = \{0\} \times \mathfrak{h}$; the converse also holds~\cite[Proposition 8.3.2]{majid2}: if a~Lie algebra
$\mathfrak{L}$ factorizes through $\mathfrak{g}$ and $\mathfrak{h}$, then there exist an isomorphism of Lie algebras
$\mathfrak{L} \cong \mathfrak{g} \bowtie \mathfrak{h}$, where $\mathfrak{g} \bowtie \mathfrak{h}$ is the bicrossed product
associated to the matched pair $(\mathfrak{g}, \mathfrak{h}, \triangleleft, \triangleright)$ whose actions are constructed from
the unique decomposition
\begin{gather}
\label{mpcan}
\left[x, g \right] = x \triangleright g + x \triangleleft g \in \mathfrak{g} + \mathfrak{h}
\end{gather}
for all $x \in \mathfrak{h}$ and $g \in \mathfrak{g}$.
The matched pair $(\mathfrak{g}, \mathfrak{h}, \triangleleft, \triangleright)$ def\/ined by~\eqref{mpcan} is called the
\emph{canonical matched pair} associated to the factorization $\mathfrak{L} = \mathfrak{g} + \mathfrak{h}$.

\begin{Remark}
%\label{complexproduct}
Over the complex numbers ${\mathbb C}$, an equivalent description for the factorization of a~Lie algebra $\mathfrak{L}$ through
two Lie subalgebras is given in~\cite[Definition~2.1]{ABDO} and in~\cite[Proposition~2.2]{AS},
in terms of \emph{complex product structures} of $\mathfrak{L}$, i.e.~linear maps $f: \mathfrak{L} \to \mathfrak{L}$
such that $f \neq \pm \Id$, $f^2 = f$ satisfying
the integrability conditions
\begin{gather*}
f ([x, y]) = [f(x), y] + [x, f(y)] - f \big([f(x), f(y)] \big)
\end{gather*}
for all $x, y \in \mathfrak{L}$.
The linear map $f: \mathfrak{g} \bowtie \mathfrak{h} \to \mathfrak{g} \bowtie \mathfrak{h}$, $f (g, h):= (g, - h)$ is a~complex
product structure on any bicrossed product $\mathfrak{g} \bowtie \mathfrak{h}$.
Conversely, if~$f$ is a~complex product structure on $\mathfrak{L}$, then $\mathfrak{L}$ factorizes through two Lie subalgebras
$\mathfrak{L} = \mathfrak{L}_{+} + \mathfrak{L}_{-}$, where $\mathfrak{L}_{\pm}$ denotes the eigenspace corresponding to the
eigenvalue $\pm 1$ of~$f$, that is $\mathfrak{L} \cong \mathfrak{L}_{+} \bowtie \mathfrak{L}_{-}$.
\end{Remark}

Let $(\mathfrak{g}, \mathfrak{h}, \triangleleft, \triangleright)$ be a~matched pair of Lie algebras.
A~linear map $r: \mathfrak{h} \to \mathfrak{g}$ is called a~\emph{deformation map}~\cite[Definition 4.1]{am-2013a} of the
matched pair $(\mathfrak{g}, \mathfrak{h}, \triangleright, \triangleleft)$ if the following compatibility holds for any $x,y\in \mathfrak{h}$
\begin{gather}
\label{factLie}
r\big([x, y]\big) - \big[r(x), r(y)\big] = r \big(y \triangleleft r(x) - x \triangleleft r(y) \big) + x \triangleright
r(y) - y \triangleright r(x).
\end{gather}
We denote by ${\mathcal D}{\mathcal M} (\mathfrak{h}, \mathfrak{g} | (\triangleright, \triangleleft))$ the set of all
deformation maps of the matched pair $(\mathfrak{g}, \mathfrak{h}, \triangleright, \triangleleft)$.
If $r \in {\mathcal D}{\mathcal M} (\mathfrak{h}, \mathfrak{g} | (\triangleright, \triangleleft))$ then $\mathfrak{h}_{r}:=
\mathfrak{h}$, as a~vector space, with the new bracket def\/ined for any~$x, y \in \mathfrak{h}$ by
\begin{gather}
\label{rLiedef}
[x, y]_{r}:= [x, y] + x \triangleleft r(y) - y \triangleleft r(x)
\end{gather}
is a~Lie algebra called the \emph{$r$-deformation} of $\mathfrak{h}$.
A~Lie algebra $\overline{\mathfrak{h}}$ is a~complement of $\mathfrak{g} \cong \mathfrak{g} \times \{0\}$ in the bicrossed
product $\mathfrak{g} \bowtie \mathfrak{h}$ if and only if $\overline{\mathfrak{h}} \cong \mathfrak{h}_{r}$, for some
deformation map $r \in {\mathcal D}{\mathcal M} (\mathfrak{h}, \mathfrak{g} | (\triangleright, \triangleleft))$ \cite[Theo\-rem~4.3]{am-2013a}.

\section{The case of perfect Lie algebras}
\label{bpbicross}

Computing all matched pairs between two given Lie algebras $\mathfrak{g}$ and $\mathfrak{h}$ and classifying all associated
bicrossed products $\mathfrak{g} \bowtie \mathfrak{h}$ is a~challenging problem.
In the case when $\mathfrak{g}:= k = k_0$, the Abelian Lie algebra of dimension~$1$, they are parameterized by
the set $\TwDer (\mathfrak{h})$ of all twisted derivations of the Lie algebra $\mathfrak{h}$
as def\/ined in~\cite[Definition~4.2]{am-2013b}: a~\emph{twisted derivation} of $\mathfrak{h}$ is a~pair $(\lambda, \Delta)$
consisting of two linear maps $\lambda: \mathfrak{h} \to k$ and $\Delta: \mathfrak{h} \to \mathfrak{h}$ such that for any~$g, h\in \mathfrak{h}$
\begin{gather}
\label{lambderivari}
\lambda ([g, h]) = 0,
\qquad
\Delta ([g, h]) = [\Delta (g), h] + [g, \Delta (h)] + \lambda(g) \Delta (h) - \lambda(h) \Delta (g).
\end{gather}
$\TwDer (\mathfrak{h})$ contains the usual space of derivations $\Der (\mathfrak{h})$ via the canonical embedding
$\Der (\mathfrak{h}) \hookrightarrow \TwDer (\mathfrak{h})$, $D \mapsto (0, D) $, which is an isomorphism if
$\mathfrak{h}$ is a~perfect Lie algebra (i.e.~$\mathfrak{h} = [\mathfrak{h}, \mathfrak{h}]$).
As a~special case of~\cite[Proposition~4.4 and Remark~4.5]{am-2013b} we have:

\begin{Proposition}
\label{mpdim1}
Let $\mathfrak{h}$ be a~Lie algebra.
Then there exists a~bijection between the set of all matched pairs $(k_0, \mathfrak{h}, \triangleleft, \triangleright)$ and the
space $\TwDer (\mathfrak{h})$ of all twisted derivations of $\mathfrak{h}$ given such that the matched pair $(k_0,
\mathfrak{h}, \triangleleft, \triangleright)$ corresponding to $(\lambda, \Delta) \in \TwDer (\mathfrak{h})$ is defined by
\begin{gather}
\label{extenddim10}
h \triangleright a= a\lambda (h),
\qquad
h \triangleleft a= a\Delta (h)
\end{gather}
for all $h \in \mathfrak{h}$ and $a \in k = k_0$.
The bicrossed product $k_0 \bowtie \mathfrak{h}$ associated to the matched pair~\eqref{extenddim10} is denoted~by
$\mathfrak{h}_{(\lambda, \Delta)}$ and has the bracket given for any~$a, b \in k$ and~$x, y \in \mathfrak{h}$ by
\begin{gather}
\{(a, x), (b, y) \}:= \big(b \lambda(x) - a\lambda(y), [x, y] + b \Delta(x) - a\Delta(y)\big).
\label{exdim300aa}
\end{gather}
A~Lie algebra $\mathfrak{L}$ contains $\mathfrak{h}$ as a~subalgebra of codimension~$1$ if and only if $\mathfrak{L}$ is
isomorphic to~$\mathfrak{h}_{(\lambda, \Delta)}$, for some $(\lambda, \Delta) \in \TwDer (\mathfrak{h})$.
\end{Proposition}

Suppose $\{e_i \,|\, i\in I \}$
is a~basis for the Lie algebra $\mathfrak{h}$.
Then, $\mathfrak{h}_{(\lambda, \Delta)}$ has $\{F, e_i \,|\, i\in I \}$ as a~basis and the bracket given for any $i\in I$~by
\begin{gather*}
%\label{extenddim200}
[e_i, F] = \lambda (e_i) F + \Delta (e_i),
\qquad
[e_i, e_j] = [e_i, e_j]_{\mathfrak{h}},
\end{gather*}
where $[-, -]_{\mathfrak{h}}$ is the bracket on $\mathfrak{h}$.
Above we identify $e_i = (0, e_i)$ and denote $F = (1, 0)$ in the bicrossed product $k_0 \bowtie \mathfrak{h}$.
Classifying the Lie algebras $\mathfrak{h}_{(\lambda, \Delta)}$ is a~dif\/f\/icult task.
In what follows we deal with this problem for a~perfect Lie algebra $\mathfrak{h}$: in this case $\TwDer (\mathfrak{h}) =
\{0\} \times \Der (\mathfrak{h})$ and we denote by $\mathfrak{h}_{(\Delta)} = \mathfrak{h}_{(0, \Delta)}$, for any $\Delta
\in \Der (\mathfrak{h})$.
\begin{Theorem}
\label{morfismeperfecte}
Let $\mathfrak{h}$ be a~perfect Lie algebra and~$\Delta, \Delta' \in \Der (\mathfrak{h})$.
Then there exists a~bijection between the set of all morphisms of Lie algebras $\varphi: \mathfrak{h}_{(\Delta)} \to
\mathfrak{h}_{(\Delta')}$ and the set of all triples $(\alpha, h, v) \in k \times \mathfrak{h} \times \Hom_{\Lie}
(\mathfrak{h}, \mathfrak{h})$ satisfying the following compatibility condition for all $x \in \mathfrak{h}$
\begin{gather}
\label{compmorfisme}
v\big(\Delta (x) \big) - \alpha \Delta ' \big(v(x) \big) = [v(x), h].
\end{gather}
The bijection is given such that the Lie algebra map $\varphi = \varphi_{(\alpha, h, v)}$ corresponding to $(\alpha, h, v)$ is
given by the formula
\begin{gather*}
%\label{compmorfisme2}
\varphi: \mathfrak{h}_{(\Delta)} \to \mathfrak{h}_{(\Delta')},
\qquad
\varphi (a, x) = (a\alpha, ah + v(x))
\end{gather*}
for all $(a, x) \in \mathfrak{h}_{(\Delta)} = k_0 \bowtie \mathfrak{h}$.
Furthermore, $\varphi = \varphi_{(\alpha, h, v)}$ is an isomorphism of Lie algebras if and only if $\alpha \neq 0$ and $v \in
\Aut_{\Lie} (\mathfrak{h})$.
\end{Theorem}

\begin{proof}
Any linear map $\varphi: k \times \mathfrak{h} \to k \times \mathfrak{h}$ is uniquely determined by a~quadruple $(\alpha, h,
\beta, v)$, where $\alpha \in k$, $h \in \mathfrak{h}$ and $\beta: \mathfrak{h} \to k$, $v: \mathfrak{h} \to \mathfrak{h}$
are~$k$-linear maps such that
\begin{gather*}
\varphi (a,x) = \varphi_{(\alpha, h, \beta, v)} = (a\alpha + \beta(x), ah + v(x)).
\end{gather*}
We will prove that~$\varphi$ def\/ined above is a~Lie algebra map if and only if~$\beta$ is the trivial map,~$v$ is a~Lie algebra
map and~\eqref{compmorfisme} holds.
It is enough to test the compatibility
\begin{gather}
\label{Liemap}
\varphi \big([(a, x), (b, y)] \big) = [\varphi(a, x), \varphi(b, y)]
\end{gather}
for all generators of $\mathfrak{h}_{(\Delta)} = k \times \mathfrak{h}$, i.e.~elements of the form $(1, 0)$ and $(0, x)$, for
all $x \in \mathfrak{h}$.
Moreover, since $\mathfrak{h}$ is perfect (i.e.~$\lambda = 0$) the bracket on $\mathfrak{h}_{(\Delta)}$ given
by~\eqref{exdim300aa} takes the form: $\{(a, x), (b, y)\}=(0, [x, y] + b \Delta(x) - a~\Delta(y))$.
Using this formula we obtain that~\eqref{Liemap} holds for $(0, x)$ and $(0, y)$ if and only if
\begin{gather*}
\beta \big([x, y]\big) = 0,
\qquad
v \big([x, y]\big) = [v(x), v(y)] + \beta (y) \Delta (v(x)) - \beta (x) \Delta (v (y)).
\end{gather*}
As $\mathfrak{h}$ is perfect these two conditions are equivalent to the fact that $\beta = 0$ and~$v$ is a~Lie algebra map.
Finally, as $\beta =0 $, we can easily show that~\eqref{Liemap} holds in $(1, 0)$ and $(0, x)$ if and only
if~\eqref{compmorfisme} holds.
Thus, we have obtained that~$\varphi$ is a~Lie algebra map if and only if~$v$ is a~Lie algebra map, $\beta = 0$
and~\eqref{compmorfisme} holds.
In what follows we denote by $\varphi_{(\alpha, h, v)}$ the Lie algebra map corresponding to a~quadruple $(\alpha, h, \beta, v)$
with $\beta = 0$.
Suppose f\/irst that $\varphi:= \varphi_{(\alpha, h, v)}$ is a~Lie algebra isomorphism.
Then, there exists a~Lie algebra map $\overline{\varphi}: = \varphi_{(\gamma, g, w)}: \mathfrak{h}_{(\Delta')} \to
\mathfrak{h}_{(\Delta)}$ such that $\varphi \circ \overline{\varphi} (a, x) = \overline{\varphi} \circ \varphi(a, x) = (a, x)$
for all $a \in k$, $x \in \mathfrak{h}$.
Thus, for all $a \in k$ and $x \in \mathfrak{h}$, we have
\begin{gather}
\label{LLL}
a\alpha \gamma = a,
\qquad
a\gamma + v(a g) + v\big(w(x)\big) = x = a\alpha g + w(a h) + w\big(v(x)\big).
\end{gather}
By the f\/irst part of~\eqref{LLL} for $a = 1$ we obtain $\alpha \gamma = 1$ and thus $\alpha \neq 0$ while the second part
of~\eqref{LLL} for $a = 0$ implies~$v$ bijective.
To end with, assume that $\alpha \neq 0$ and $v \in \Aut_{\Lie} (\mathfrak{h})$.
Then, it is straightforward to see that $\varphi = \varphi_{(\alpha, h, v)}$ is an isomorphism with the inverse given~by
$\varphi^{-1}: = \varphi_{(\alpha^{-1}, -\alpha^{-1} v^{-1}(h), v^{-1})}$.
\end{proof}

Let $k^*$ be the units group of~$k$ and $(\mathfrak{h}, +)$ the underlying Abelian group of the Lie algebra $\mathfrak{h}$.
Then the map given for any $\alpha \in k^*$, $v\in \Aut_{\Lie} (\mathfrak{h})$ and $h\in \mathfrak{h}$ by
\begin{gather*}
\varphi: \ k^* \times \Aut_{\Lie} (\mathfrak{h}) \to \Aut_{\Gr} (\mathfrak{h}, +),
\qquad
\varphi (\alpha, v) (h):= \alpha^{-1} v(h)
\end{gather*}
is a~morphism of groups.
Thus, we can construct the semidirect product of groups $\mathfrak{h} \ltimes_{\varphi}(k^* \times \Aut_{\Lie}(\mathfrak{h}))$ associated to~$\varphi$.
The next result shows that $\Aut_{\Lie} (\mathfrak{h}_{(\Delta)})$ is isomorphic to a~certain subgroup of the semidirect
product of groups $\mathfrak{h} \ltimes_{\varphi}(k^* \times \Aut_{\Lie} (\mathfrak{h}))$.

\begin{Corollary}
\label{izoaut}
Let $\mathfrak{h}$ be a~perfect Lie algebra and~$\Delta, \Delta' \in \Der (\mathfrak{h})$.
Then the Lie algebras $\mathfrak{h}_{(\Delta)}$ and $\mathfrak{h}_{(\Delta')}$ are isomorphic if and only if there exists
a~triple $(\alpha, h, v) \in k^* \times \mathfrak{h} \times \Aut_{\Lie} (\mathfrak{h})$ such that $v \circ \Delta -
\alpha \Delta ' \circ v = [v(-), h]$.
Furthermore, there exists an isomorphism of groups
\begin{gather*}
\Aut_{\Lie} (\mathfrak{h}_{(\Delta)}) \cong {\mathcal G} (\mathfrak{h}, \Delta):= \{(\alpha, h, v) \in k^* \times
\mathfrak{h} \times \Aut_{\Lie} (\mathfrak{h}) \,|\, v \circ \Delta - \alpha \Delta \circ v = [v(-), h] \},
\end{gather*}
where ${\mathcal G} (\mathfrak{h}, \Delta)$ is a~group with respect to the following multiplication
\begin{gather}
\label{graut}
(\alpha, h, v) \cdot (\beta, g, w):= (\alpha \beta, \beta h + v(g), v \circ w)
\end{gather}
for all $(\alpha, h, v), (\beta, g, w)  \in {\mathcal G} (\mathfrak{h}, \Delta)$.
Moreover, the canonical map
\begin{gather*}
{\mathcal G} (\mathfrak{h}, \Delta) \longrightarrow \mathfrak{h} \ltimes_{\varphi} \big(k^* \times \Aut_{\Lie}
(\mathfrak{h}) \big),
\qquad
(\alpha, h, v) \mapsto \big (\alpha^{-1} h, (\alpha, v) \big)
\end{gather*}
is an injective morphism of groups.
\end{Corollary}

\begin{proof}
The f\/irst part follows trivially from Theorem~\ref{morfismeperfecte}.
Consider now~$\gamma, \psi \in \Aut_{\Lie} (\mathfrak{h}_{(\Delta)})$.
Using again Theorem~\ref{morfismeperfecte}, we can f\/ind $(\alpha, h, v),
(\beta, g, w) \in k^* \times \mathfrak{h} \times \Aut_{\Lie} (\mathfrak{h})$
such that $\gamma = \varphi_{(\alpha, h, v)}$ and $\psi = \varphi_{(\beta, g, w)}$.
Then, for all $a \in k$, $x \in \mathfrak{h}$ we have
\begin{gather*}
\varphi_{(\alpha, h, v)} \circ \varphi_{(\beta, g, w)} (a, x)  =  \varphi_{(\alpha, h, v)} \big(a \beta, ag + w(x)\big)
 =  \big(\alpha \beta a, a\beta h + av(g) + v \circ w(x)\big)
\\
\hphantom{\varphi_{(\alpha, h, v)} \circ \varphi_{(\beta, g, w)} (a, x)}{}
 =  \varphi_{(\alpha \beta, \beta h + v(g), v \circ w)}(a, x).
\end{gather*}
Thus, $\Aut_{\Lie} (\mathfrak{h}_{(\Delta)}) $ is isomorphic to ${\mathcal G} (\mathfrak{h}, \Delta)$ with the
multiplication given by~\eqref{graut}.
The last assertion follows by a~routine computation.
\end{proof}

\begin{Remark}
\label{cazinner}
Let $\Delta = [x_0, -]$ be an inner derivation of a~perfect Lie algebra $\mathfrak{h}$.
Then the group $\Aut_{\Lie} (\mathfrak{h}_{([x_0, -])})$ admits a~simpler description as follows
\begin{gather*}
{\mathcal G} (\mathfrak{h}, [x_0, -]) = \{(\alpha, h, v) \in k^* \times \mathfrak{h} \times \Aut_{\Lie} (\mathfrak{h}) \,|\,
v(x_0) - \alpha x_0 + h \in {\rm Z} (\mathfrak{h})\},
\end{gather*}
where ${\rm Z} (\mathfrak{h})$ is the center of $\mathfrak{h}$.
Assume in addition that $\mathfrak{h}$ has trivial center, i.e.~${\rm Z} (\mathfrak{h}) = \{0\}$; it follows that there exists
an isomorphism of groups
\begin{gather*}
\Aut_{\Lie} (\mathfrak{h}_{([x_0, -])}) \cong k^* \times \Aut_{\Lie} (\mathfrak{h}),
\end{gather*}
since in this case any element~$h$ from a~triple $(\alpha, h, v) \in {\mathcal G} (\mathfrak{h}, [x_0, -])$ must be equal to $
\alpha x_0 - v (x_0)$.
Moreover, in this context, the multiplication given by~\eqref{graut} is precisely that of a~direct product of groups.
\end{Remark}

A Lie algebra $\mathfrak{h}$ is called \emph{complete} (see~\cite{JMZ, SZ} for examples and structural results on this class of
Lie algebras) if $\mathfrak{h}$ has trivial center and any derivation is inner.
A~complete and perfect Lie algebra is called \emph{sympathetic}~\cite{ben}: semisimple Lie algebras over a~f\/ield of
characteristic zero are sympathetic and there exists a~sympathetic non-semisimple Lie algebra in dimension $25$.
For sympathetic Lie algebras, Theorem~\ref{morfismeperfecte} takes the following form which considerably
improves~\cite[Corollary~4.10]{am-2013b}, where the classif\/ication is made only up to an isomorphism of Lie algebras which acts
as identity on~$\mathfrak{h}$.

\begin{Corollary}
\label{morfismeperfectea}
Let $\mathfrak{h}$ be a~sympathetic Lie algebra.
Then up to an isomorphism of Lie algebras there exists only one Lie algebra that contains $\mathfrak{h}$ as a~Lie subalgebra of
codimension one, namely the direct product $k_0 \times \mathfrak{h}$ of Lie algebras.
Furthermore, there exists an isomorphism of groups $\Aut_{\Lie} (k_0 \times \mathfrak{h}) \cong k^* \times \Aut_{\Lie} (\mathfrak{h})$.
\end{Corollary}

\begin{proof}
Since $\mathfrak{h}$ is perfect any Lie algebra that contains $\mathfrak{h}$ as a~Lie subalgebra of codimension $1$ is
isomorphic to $\mathfrak{h}_{(D)}$, for some $D \in \Der (\mathfrak{h})$.
As $\mathfrak{h}$ is also complete, any derivation is inner.
For an arbitrary derivation $D = [d, -]$ we can prove that $\mathfrak{h}_{(D)} \cong \mathfrak{h}_{(0)}$, where $0 = [0, -]$ is
the trivial derivation and moreover $\mathfrak{h}_{(0)}$ is just the direct product of Lie algebras $k_0 \times \mathfrak{h}$.
Indeed, by taking $(\alpha, h, v):= (1, -d, \Id_{\mathfrak{h}})$ one can see that relation~\eqref{compmorfisme} holds for
$D = [d, -]$ and $D' = [0, -]$, that is $\mathfrak{h}_{(D)} \cong \mathfrak{h}_{(0)}$.
The f\/inal part follows from Remark~\ref{cazinner}.
\end{proof}

\begin{Remark}
\label{defhdelta}
Let $\mathfrak{h}$ be a~perfect Lie algebra with a~basis $\{e_i \,|\, i\in I \}$, $\Delta \in \Der (\mathfrak{h})$ a~given
derivation and consider the extension $k_0 \subseteq \mathfrak{h}_{(\Delta)} = k_0 \bowtie \mathfrak{h}_{(\Delta)}$.
In order to determine all complements of $k_0$ in $\mathfrak{h}_{(\Delta)}$ we have to describe the set of all deformation maps
$r: \mathfrak{h} \to k_0$ of the matched pair~\eqref{extenddim10}.
A~deformation map is completely determined by a~family of scalars $(a)_{i\in I}$ satisfying the following compatibility
condition for any~$i, j\in I$
\begin{gather*}
r \big([e_i, e_j]_{\mathfrak{h}} \big) = r \big(a_i \Delta (e_j) - a_j \Delta (e_i) \big)
\end{gather*}
via the relation $r (e_i) = a_i$.
For such an $r = (a_i)_{i\in I}$, the~$r$-deformation of $\mathfrak{h}$ is the Lie algebra $\mathfrak{h}_r$ having $\{e_i \,|\, i\in I\}$
as a~basis and the bracket def\/ined for any~$i, j\in I$ by
\begin{gather*}
[e_i, e_j]_r = [e_i, e_j]_{\mathfrak{h}} + a_j \Delta (e_i) - a_i \Delta (e_j).
\end{gather*}
Any complement of $k_0$ in $\mathfrak{h}_{(\Delta)}$ is isomorphic to such an $\mathfrak{h}_r$.
An explicit example in dimension $5$ is given below.
\end{Remark}

\begin{Example}
\label{prefectnecomplet}
Let~$k$ be a~f\/ield of characteristic $\neq 2$ and $\mathfrak{h}$ the perfect $5$-dimensional Lie algebra with a~basis $\{e_{1},
e_{2}, e_{3}, e_{4}, e_{5}\}$ and bracket given by
\begin{alignat*}{4}
& [e_{1}, e_{2}] = e_{3},
\qquad &&
[e_{1}, e_{3}] = -2e_{1},
\qquad &&
[e_{1}, e_{5}] = [e_{3}, e_{4}] = e_{4},&
\\
&  [e_{2}, e_{3}  ] = 2e_{2},
\qquad &&
 [e_{2}, e_{4}  ] = e_{5},
\qquad &&
 [e_{3}, e_{5}  ] = - e_{5}.&
\end{alignat*}
By a~straightforward computation it can be proved that the space of derivations $\Der(\mathfrak{h})$ coincides with the
space of all matrices from $\mathcal{M}_{5}(k)$ of the form
\begin{gather*}
A= \left(
\begin{matrix} a_{1} & 0 & -2a_4 & 0 & 0
\\
0 & -a_{1} & -2a_{2} & 0 & 0
\\
a_{2} & a_{4} & 0 & 0 & 0
\\
a_{3} & 0 & a_{5} & a_{6} & a_{4}
\\
0 & a_{5} & -a_{3} & -a_{2} & (a_{6}-a_{1})
\end{matrix}
\right)
\end{gather*}
for all $a_1, \dots, a_6 \in k$.
Thus $\mathfrak{h}$ is not complete since $\Der(\mathfrak{h})$ has dimension~$6$.
One can show easily that the derivation $\Delta:= e_{11} - e_{41} - e_{22} + e_{53} - e_{44} - 2 e_{55}$ is not inner, where
$e_{i j} \in \mathcal{M}_{n}(k)$ is the matrix having $1$ in the $(i,j)^{\rm th}$ position and zeros elsewhere.
For the derivation~$\Delta$ we consider the extension $k_0 \subseteq k_0 \bowtie \mathfrak{h} = \mathfrak{h}_{(\Delta)}$ and we
will describe all the complements of $k_0$ in $\mathfrak{h}_{(\Delta)}$.
By a~routine computation it can be seen that $r: \mathfrak{h} \to k_0$ is a~deformation map of the matched
pair~\eqref{extenddim10} if and only if $r:= 0$ (the trivial map) or~$r$ is given~by
\begin{gather*}
r(e_1):= a,
\qquad
r(e_2):= - a^{-1},
\qquad
r(e_3) = 2,
\qquad
r(e_4) = r(e_5) = 0
\end{gather*}
for some $a \in k^*$.
Thus a~Lie algebra $\mathfrak{C}$ is a~complement of $k_0$ in $\mathfrak{h}_{(\Delta)}$ if and only if $\mathfrak{C} \cong
\mathfrak{h}$ or $\mathfrak{C} \cong \mathfrak{h}_a$, where $\mathfrak{h}_a$ is the $5$-dimensional Lie algebra with basis
$\{e_{1}, e_{2}, e_{3}, e_{4}, e_{5}\}$ and bracket given~by
\begin{gather*}
 [e_1, e_2  ]_a:= - a^{-1} e_1 + ae_2 + e_3 + a^{-1} e_4,
\qquad
 [e_1, e_3  ]_a:= -2 e_4 - ae_5,
\qquad
 [e_1, e_4  ]_a:= ae_4,
\\
 [e_1, e_5  ]_a:= e_4 + 2 ae_5,
\qquad
 [e_2, e_3  ]_a:= a^{-1} e_5,
\qquad
 [e_2, e_4  ]_a:= e_5 - a^{-1} e_4,
\\
 [e_2, e_5 ]_a:= -2 a^{-1} e_5,
\qquad
 [e_3, e_4  ]_a:= 3 e_4,
\qquad
 [e_3, e_5  ]_a:= 3 e_5
\end{gather*}
for any $a \in k^*$.
Remark that none of the matched pair deformations $\mathfrak{h}_a$ of the Lie algebra $\mathfrak{h}$ is perfect since the
dimension of the derived algebra $[\mathfrak{h}_a, \mathfrak{h}_a]$ is equal to~$3$.
\end{Example}

\section{The non-perfect case}
\label{mpdefo}

In Section~\ref{bpbicross} we have described and classif\/ied all bicrossed products $k_0 \bowtie \mathfrak{h}$ for a~perfect Lie
algebra $\mathfrak{h}$; furthermore, Remark~\ref{defhdelta} and Example~\ref{prefectnecomplet} describe all complements of $k_0$
in a~given bicrossed product $k_0 \bowtie \mathfrak{h}$.
In this section we approach the same questions for a~given non-perfect Lie algebra $\mathfrak{h}:= \mathfrak{l} (2n+1, k)$,
where $\mathfrak{l} (2n+1, k)$ is the $(2n+1)$-dimensional Lie algebra with basis $\{E_i, F_i, G \,|\, i = 1, \dots, n\}$ and
bracket given for any $i = 1, \dots, n$~by
\begin{gather*}
[E_i, G]:= E_i,
\qquad
[G, F_i]:= F_i.
\end{gather*}
First, we shall describe all bicrossed products $k_0 \bowtie \mathfrak{l} (2n+1, k)$: they will explicitly describe all Lie
algebras which contain $\mathfrak{l} (2n+1, k)$ as a~subalgebra of codimension~$1$.
Then, as the second step, we shall f\/ind all~$r$-deformations of the Lie algebra $\mathfrak{l} (2n+1, k)$, for two given
extensions $k_0 \subseteq k_0 \bowtie \mathfrak{l} (2n+1, k)$.
Based on Proposition~\ref{mpdim1} we have to compute f\/irst the space $\TwDer(\mathfrak{l} (2n+1, k))$ of all
twisted derivations.

\begin{Proposition}
\label{tw2n1}
There exists a~bijection between $\TwDer(\mathfrak{l} (2n+1, k))$ and the set of all matrices $(A, B, C, D,
\lambda_0, \delta) \in {\rmM}_n (k)^4 \times k \times k^{2n+1}$ satisfying the following conditions
\begin{gather}
\lambda_0 A=- \delta_{2n+1} I_n,
\qquad
(2 + \lambda_0) B = 0,
\qquad
(2 - \lambda_0) C = 0,
\qquad
\lambda_0 D = \delta_{2n+1} I_n,
\label{primaa}
\end{gather}
where $\delta = (\delta_1, \dots, \delta_{2n+1}) \in k^{2n+1}$.
The bijection is given such that the twisted derivation $(\lambda, \Delta) \in \TwDer (\mathfrak{l} (2n+1, k))$
associated to $(A, B, C, D, \lambda_0, \delta)$ is given by
\begin{gather}
\lambda (E_i) = \lambda (F_i):= 0,
\qquad
\lambda (G):= \lambda_0,
\label{primab}
\\
\Delta:=
\begin{pmatrix}
A & B & \delta_1
\\
C & D  :
\\
0 & 0 & \delta_{2n+1}
\end{pmatrix}.
\label{primac}
\end{gather}
${\mathcal T} (n)$ denotes the set of all $(A, B, C, D, \lambda_0, \delta) \in {\rmM}_n (k)^4 \times k \times k^{2n+1}$
satisfying~\eqref{primaa}.
\end{Proposition}

\begin{proof}
The f\/irst compatibility condition~\eqref{lambderivari} shows that a~linear map $\lambda: \mathfrak{l} (2n+1, k) \to k$ of
a~twisted derivation $(\lambda, D)$ must have the form given by~\eqref{primab}, for some $\lambda_0 \in k$.
We shall f\/ix such a~map for a~given $\lambda_0 \in k$.
We write down the linear map $\Delta: \mathfrak{l} (2n+1, k) \to \mathfrak{l} (2n+1, k)$ as a~matrix associated to the basis
$\{E_1, \dots, E_n, F_1, \dots, F_n, G \}$ of $\mathfrak{l} (2n+1, k)$, as follows
\begin{gather*}
\Delta =
\begin{pmatrix}
A & B & d_{1, 2n+1}
\\
C & D  :
\\
d_{2n+1, 1}  ..
  d_{2n+1, 2n+1}
\end{pmatrix}
\end{gather*}
for some matrices~$A,B,C, D \in {\rmM}_n (k)$ and some scalars $d_{i, j} \in k$, for all~$i, j = 1, \dots, 2n+1$.
We denote $A = (a_{ij})$, $B = (b_{ij})$, $C = (c_{ij})$, $D = (d_{ij})$.
It remains to check the compatibility condition~\eqref{lambderivari} for~$\Delta$, i.e.\ %the following compatibility
\begin{gather*}
\Delta ([g, h]) = [\Delta(g), h] + [g, \Delta (h)] + \lambda(g) \Delta (h) - \lambda(h) \Delta (g)
\end{gather*}
for all $g \neq h \in \{E_1, \dots, E_n, F_1, \dots, F_n, G\}$.
As this is a~routinely straightforward computation we will only indicate the main steps of the proof.
We can easily see that the compatibility condition~\eqref{lambderivari} holds
for $(g, h) = (E_i, E_j)$ if and only if $d_{2n+1,i} = 0$, for all $i = 1, \dots, n$.
In the same way~\eqref{lambderivari} holds for $(g, h) = (F_i, F_j)$ if and only if $d_{2n+1, n + i} = 0$, for all $i = 1, \dots, n$.
This shows that $\Delta $ has the form~\eqref{primac}, that is the f\/irst $2n$ entries from the last row of the matrix~$\Delta$
are all zeros and we will denote the last column of~$D$ by $(d_{1, 2n+1}, \dots, d_{2n+1, 2n+1}) = \delta = (\delta_1, \dots,
\delta_{2n+1})$.
It follows from here that~\eqref{lambderivari} holds trivially for the pair $(g, h) = (E_i, F_j)$.
An easy computation shows that~\eqref{lambderivari} holds for $(g, h) = (E_i, G)$ if and only if the following equation holds
\begin{gather*}
(1 - \lambda_0) \left(\sum\limits_{j=1}^n a_{j, i} E_j + \sum\limits_{j=1}^n c_{j, i} F_j \right) = \sum\limits_{j=1}^n a_{j, i}
E_j - \sum\limits_{j=1}^n c_{j, i} F_j + \delta_{2n+1} E_i,
\end{gather*}
which is equivalent to $- \lambda_0 A= \delta_{2n+1} I_n$ and $(2 - \lambda_0) C = 0$, i.e.~the f\/irst and the third equations
from~\eqref{primaa}.
A~similar computation shows that~\eqref{lambderivari} holds for $(g, h) = (G, F_i)$ if and only if $(2 + \lambda_0) B = 0$ and
$\lambda_0 D = \delta_{2n+1} I_n$ and the proof is f\/inished.
\end{proof}

Let $\mathfrak{l} (2n+1, k)_{(A, B, C, D, \lambda_0, \delta)}$ be the bicrossed product $k_0 \bowtie \mathfrak{l} (2n+1, k)$
associated to the matched pair given by
the twisted derivation $\big(A = (a_{ji}), B = (b_{ji}), C = (c_{ji}), D = (d_{ji}),\lambda_0$, $\delta = (\delta_j) \big) \in {\mathcal T} (n)$.
From now on we will use the following convention: if one of the elements of the $6$-tuple ($A$,~$B$,~$C$,~$D$,
$\lambda_0$,~$\delta$) is equal to $0$ then we will omit it when writing down the Lie algebra $\mathfrak{l} (2n+1, k)_{(A,B,C,D,\lambda_0,\delta)}$.
A~basis of $\mathfrak{l} (2n+1, k)_{(A, B, C, D, \lambda_0, \delta)}$ will be denoted by $\{E_i, F_i, G, H \,|\, i = 1, \dots, n\}$:
these Lie algebras can be explicitly described by f\/irst computing the set ${\mathcal T} (n)$ and then using
Proposition~\ref{mpdim1}.
Considering the equations~\eqref{primaa} which def\/ine ${\mathcal T}(n)$ a~discussion involving the f\/ield~$k$ and the scalar~$\lambda_0$ is mandatory.
For two sets~$X$ and~$Y$ we shall denote by $X \sqcup Y$ the disjoint union of~$X$ and~$Y$.
As a~conclusion of the above results we obtain:
\begin{Theorem}
\label{teorema11}
$(1)$ If~$k$ is a~field such that $\charop (k) \neq 2$ then
\begin{gather*}
{\mathcal T} (n) \cong \big((k \setminus \{0, \pm 2\}) \times k^{2n+1} \big) \sqcup \big ({\rmM}_n(k)^2 \times k^{2n} \big)
\sqcup \big({\rmM}_n(k) \times k^{2n+1} \big) \sqcup \big({\rmM}_n(k) \times k^{2n+1} \big)
\end{gather*}
and the four families of Lie algebras containing $\mathfrak{l} (2n+1, k)$ as a~subalgebra of codimension $1$ are the following:

$\bullet$ the Lie algebra $\mathfrak{l}^1 (2n+1, k)_{(\lambda_0, \delta)}$ with the bracket given for any $i = 1, \dots, n$ by
\begin{gather*}
 [E_i, G  ] = E_i,
\qquad
 [G, F_i  ] = F_i,
\qquad
 [E_i, H  ] = - \lambda_0^{-1} \delta_{2n+1} E_i,
\\
 [F_i, H  ] = \lambda_0^{-1} \delta_{2n+1} F_i,
\qquad
 [G, H  ] = \lambda_0 H + \sum\limits_{j=1}^n \delta_j E_j + \sum\limits_{j=1}^n \delta_{n +j} F_j + \delta_{2n+1} G
\end{gather*}
for all $(\lambda_0, \delta) \in (k \setminus \{0, \pm 2\}) \times k^{2n+1}$.

$\bullet$ the Lie algebra $\mathfrak{l}^2 (2n+1, k)_{(A, D, \delta)}$ with the bracket given for any $i = 1, \dots, n$ by
\begin{gather*}
 [E_i, G  ] = E_i,
\qquad
 [G, F_i  ] = F_i,
\qquad
 [E_i, H  ] = \sum\limits_{j=1}^n a_{ji} E_j,
\\
 [F_i, H  ] = \sum\limits_{j=1}^n d_{ji} F_j,
\qquad
 [G, H  ] = \sum\limits_{j=1}^n \delta_j E_j + \sum\limits_{j=1}^n \delta_{n +j} F_j
\end{gather*}
for all $(A = (a_{ij}), D = (d_{ij}), \delta) \in {\rmM}_n(k) \times {\rmM}_n(k) \times k^{2n}$.

$\bullet$ the Lie algebra $\mathfrak{l}^3 (2n+1, k)_{(C, \delta)}$ with the bracket given for any $i = 1, \dots, n$ by
\begin{gather*}
 [E_i, G  ] = E_i,
\qquad
 [G, F_i  ] = F_i,
\qquad
 [E_i, H  ] = - 2^{-1} \delta_{2n+1} E_i + \sum\limits_{j=1}^n c_{ji} F_j,
\\
 [F_i, H  ] = 2^{-1} \delta_{2n+1} F_i,
\qquad
 [G, H  ] = 2 H + \sum\limits_{j=1}^n \delta_j E_j + \sum\limits_{j=1}^n \delta_{n +j} F_j + \delta_{2n+1} G
\end{gather*}
for all $(C = (c_{ij}), \delta) \in {\rmM}_n(k) \times k^{2n+1}$.

$\bullet$ the Lie algebra $\mathfrak{l}^4 (2n+1, k)_{(B, \delta)} $ with the bracket given for any $i = 1, \dots, n$ by
\begin{gather*}
 [E_i, G  ] = E_i,
\qquad
 [G, F_i  ] = F_i,
\qquad
 [F_i, H  ] = \sum\limits_{j=1}^n b_{ji} E_j - 2^{-1} \delta_{2n+1} F_i,
\\
 [E_i, H ] = 2^{-1} \delta_{2n+1} E_i,
\qquad
 [G, H  ] = - 2 H + \sum\limits_{j=1}^n \delta_j E_j + \sum\limits_{j=1}^n \delta_{n +j} F_j + \delta_{2n+1} G
%\label{ba44}
\end{gather*}
for all $(B = (b_{ij}), \delta) \in {\rmM}_n(k) \times k^{2n+1}$.

$(2)$ If $\charop (k) = 2$ then
\begin{gather*}
{\mathcal T} (n) \cong \big({\rmM}_n(k)^4 \times k^{2n} \big) \sqcup \big (k^* \times k^{2n+1} \big)
\end{gather*}
and the two families of Lie algebras containing $\mathfrak{l} (2n+1, k)$ as a~subalgebra of codimension $1$ are the following:

$\bullet$ the Lie algebra $\mathfrak{l}_1 (2n+1, k)_{(A, B, C, D, \delta)}$ with the bracket given for any $i = 1, \dots, n$ by
\begin{gather*}
 [E_i, G  ] = E_i,
\qquad
 [G, F_i  ] = F_i,
\qquad
 [E_i, H  ] = \sum\limits_{j = 1}^n \big(a_{ji} E_j + c_{ji} F_j \big),
\\
 [F_i, H  ] = \sum\limits_{j = 1}^n \big(b_{ji} E_j + d_{ji} F_j \big),
\qquad
 [G, H  ] = \sum\limits_{j=1}^n \delta_j E_j + \sum\limits_{j=1}^n \delta_{n +j} F_j
\end{gather*}
for all $(A, B, C, D, \delta) \in {\rmM}_n(k)^4 \times k^{2n}$.

$\bullet$ the Lie algebra $\mathfrak{l}_2 (2n+1, k)_{(\lambda_0, \delta)}$ with the bracket given for any $i = 1, \dots, n$ by
\begin{gather*}
 [E_i, G  ] = E_i,
\qquad
 [G, F_i  ] = F_i,
\qquad
 [E_i, H  ] = - \lambda_0^{-1} \delta_{2n+1} E_i,
\\
 [F_i, H  ] = \lambda_0^{-1} \delta_{2n+1} F_i,
\qquad
 [G, H  ] = \lambda_0 H + \sum\limits_{j=1}^n \delta_j E_j + \sum\limits_{j=1}^n \delta_{n +j} F_j + \delta_{2n+1} G
%\label{ba441}
\end{gather*}
for all $(\lambda_0, \delta) \in k^* \times k^{2n+1}$.
\end{Theorem}

\begin{proof}
The proof relies on the use of Propositions~\ref{mpdim1} and~\ref{tw2n1} as well as the equations~\eqref{primaa}
def\/ining ${\mathcal T} (n)$.
Besides the discussion on the characteristic of~$k$ it is also necessary to consider whether $\lambda_0$ belongs to the set
$\{0, 2, -2\}$.
In the case that $\charop (k) \neq 2$, the f\/irst Lie algebra listed is the bicrossed product which corresponds to the case
when $\lambda_0 \notin \{0, 2, -2\}$.
In this case, we can easily see that $\big(A, B, C, D, \lambda_0, \delta = (\delta_j) \big) \in {\mathcal T} (n)$ if and only
if $B = C = 0$, $A = - \lambda_0^{-1} \delta_{2n+1} I_n$ and $D = \lambda_0^{-1} \delta_{2n+1} I_n$.
The Lie algebra $\mathfrak{l}^1 (2n+1, k)_{(\lambda_0, \delta)}$ is exactly the bicrossed product $k_0 \bowtie \mathfrak{l}
(2n+1, k)$ corresponding to this twisted derivation.
The Lie algebra $\mathfrak{l}^2 (2n+1, k)_{(A, D, \delta)}$ is the bicrossed product $k_0 \bowtie \mathfrak{l} (2n+1, k)$
corresponding to the case $\lambda_0 = 0$ while the last two Lie algebras are the bicrossed products $k_0 \bowtie \mathfrak{l}
(2n+1, k)$ associated to the case when $\lambda_0 = 2$ and respectively $\lambda_0 = -2$.

If the characteristic of~$k$ is equal to $2$ we distinguish the following two possibilities: the Lie algebra $\mathfrak{l}_1
(2n+1, k)_{(A, B, C, D, \delta)}$ is the bicrossed product $k_0 \bowtie \mathfrak{l} (2n+1, k)$ associated to $\lambda_0 = 0$
while the Lie algebra $\mathfrak{l}_2 (2n+1, k)_{(\lambda_0, \delta)}$ is the same bicrossed product but associated to
$\lambda_0 \neq 0$.
\end{proof}

Let~$k$ be a~f\/ield of characteristic $\neq 2$ and $\mathfrak{l}^1 (2n+1, k)_{(\lambda_0, \delta)}$ the Lie algebra of
Theorem~\ref{teorema11}.
In order to keep the computations ef\/f\/icient we will consider $\lambda_0:= 1$ and $\delta:= (0, \dots, 0, 1)$ and we denote by $L
(2n+2, k):= \mathfrak{l}^1 (2n+1, k)_{(1, (0, \dots, 0, 1))}$, the $(2n+2)$-dimensional Lie algebra having a~basis $\{E_i, F_i,
G, H \,|\, i = 1, \dots, n \}$ and the bracket def\/ined for any $i = 1, \dots, n$ by
\begin{gather*}
[E_i, G ] = E_i,
\qquad
[G, F_i] = F_i,
\qquad
[E_i, H ] = - E_i,
\qquad
[F_i, H ] = F_i,
\qquad
[G, H ] = H + G.
%\label{cazulspec}
\end{gather*}

We consider the Lie algebra extension $kH \subset L (2n+2, k)$, where $kH \cong k_0$ is the Abelian Lie algebra of dimension
$1$.
Of course, $L (2n+2, k)$ factorizes through $kH$ and $\mathfrak{l} (2n+1, k)$,
i.e.~$L (2n+2, k) = kH \bowtie \mathfrak{l} (2n+1,k)$~-- the actions $\triangleleft: \mathfrak{l} (2n+1, k) \times kH \to \mathfrak{l} (2n+1, k)$
and $\triangleright: \mathfrak{l}(2n+1, k) \times kH \to kH$ of the canonical matched pair are given by
\begin{gather}
E_i \triangleleft H:= - E_i,
\qquad
F_i \triangleleft H:= F_i,
\qquad
G \triangleleft H:= G,
\qquad
G\triangleright H:= H
\label{mpcanon}
\end{gather}
and all undef\/ined actions are zero.
Next we compute the set ${\mathcal D}{\mathcal M}(\mathfrak{l} (2n+1, k), kH | (\triangleright, \triangleleft))$ of
all deformation maps of the matched pair $(kH, \mathfrak{l} (2n+1, k), \triangleright, \triangleleft)$ given by~\eqref{mpcanon}.

\begin{Lemma}
\label{defmaps}
Let~$k$ be a~field of characteristic $\neq 2$.
Then there exists a~bijection
\begin{gather*}
{\mathcal D}{\mathcal M} \big(\mathfrak{l} (2n+1, k), kH | (\triangleright, \triangleleft) \big) \cong \big(k^n \setminus
\{0\} \big) \sqcup \big(k^n \times k \big).
\end{gather*}
The bijection is given such that the deformation map $r = r_a: \mathfrak{l} (2n+1, k) \to kH$ associated to $a = (a_i) \in k^n
\setminus \{0\}$ is given~by
\begin{gather}
\label{def1a}
r (E_i):= a_i H,
\qquad
r(F_i):= 0,
\qquad
r (G):= H,
\end{gather}
while the deformation map $r = r_{(b, c)}: \mathfrak{l} (2n+1, k) \to kH$ associated to $(b = (b_i), c) \in k^n \times k$ is
given as follows
\begin{gather}
\label{def2a}
r (E_i):= 0,
\qquad
r(F_i):= b_i H,
\qquad
r (G):= c H
\end{gather}
for all $i = 1, \dots, n$.
\end{Lemma}
\begin{proof}
Any linear map $r: \mathfrak{l} (2n+1, k) \to kH$ is uniquely determined by a~triple $(a = (a_i), b = (b_i), c) \in k^n \times
k^n \times k$ via: $r(E_i):= a_i H$, $r(F_i):= b_i H$ and $r(G):= c H$, for all $i = 1, \dots, n$.
We need to check under what conditions such a~map $r = r_{(a, b, c)}$ is a~deformation map.
Since $kH$ is Abelian, equation~\eqref{factLie} comes down to
\begin{gather}
\label{defabc}
r([x, y]) = r \big(y \triangleleft r(x) - x \triangleleft r(y) \big) + x \triangleright r(y) - y \triangleright r(x),
\end{gather}
which needs to be checked for all~$x, y \in \{E_i, F_i, G \,|\, i = 1, \dots, n \}$.
Notice that~\eqref{defabc} is symmetrical i.e.~if~\eqref{defabc} is fulf\/illed for $(x, y)$ then~\eqref{defabc} is also fulf\/illed
for $(y, x)$.
By a~routinely computation it can be seen that $r = r_{(a, b, c)}$ is a~deformation map if and only if
\begin{gather}
\label{defabc11}
a_i b_j = 0,
\qquad
(1 - c) a_i = 0
\end{gather}
for all~$i, j = 1, \dots, n$.
Indeed,~\eqref{defabc} holds for $(x, y) = (E_i, F_j)$ if and only if $a_i b_j = 0$ and it holds for $(x, y) = (E_i, G)$ if and
only if $a_i = a_i c$.
The other cases left to study are either automatically fulf\/illed or equivalent to one of the two conditions above.
The f\/irst condition of~\eqref{defabc11} divides the description of deformation maps into two cases: the f\/irst one corresponds to
$a = (a_i) \neq 0$ and we automatically have $b = 0$ and $c = 1$.
The second case corresponds to $a:= 0$ which implies that~\eqref{defabc11} holds for any $(b, c) \in k^n \times k$.
\end{proof}

The next result describes all deformations of $\mathfrak{l} (2n+1, k)$ associated to the canonical matched pair $(kH,
\mathfrak{l} (2n+1, k), \triangleright, \triangleleft)$ given by~\eqref{mpcanon}.

\begin{Proposition}
\label{mpdef}
Let~$k$ be a~field of characteristic $\neq 2$ and the extension of Lie algebras $kH \subset L (2n+2, k)$.
Then a~Lie algebra $\mathfrak{C}$ is a~complement of $kH$ in $L (2n+2, k)$ if and only if $\mathfrak{C}$ is isomorphic to one of
the Lie algebras from the three families defined below:

$\bullet$ the Lie algebra $\mathfrak{l}_{(a)} (2n+1, k)$ having the bracket def\/ined for any $i = 1, \dots, n$~by{\samepage
\begin{gather}
\label{1famdef}
[E_i, E_j]_a:= a_i E_j - a_j E_i,
\qquad
[E_i, F_j]_a:= -a_i F_j,
\qquad
[E_i, G]_a:= - a_i G
\end{gather}
for all $a = (a_i) \in k^n \setminus \{0\}$.}

$\bullet$ the Lie algebra $\mathfrak{l}'_{(b)} (2n+1, k)$ having the bracket def\/ined for any $i = 1, \dots, n$ by
\begin{gather*}
[E_i, F_j]_{b}:= - b_j E_i,
\qquad
[E_i, G]_{b}:= - E_i,
%\label{2afamdef}
\\
[F_i, F_j]_{b}:= b_j F_i - b_i F_j,
\qquad
[F_i, G]_{b}:= F_i - b_i G
%\label{2bfamdef}
\end{gather*}
for all $b = (b_i) \in k^n$.

$\bullet$ the Lie algebra $\mathfrak{l}''_{(b)} (2n+1, k)$ having the bracket def\/ined for any $i = 1, \dots, n$~by
\begin{gather*}
[E_i, F_j]_{b}:= - b_j E_i,
\qquad
[F_i, F_j]_{b}:= b_j F_i - b_i F_j,
\qquad
[F_i, G]_{b}:= - b_i G
%\label{2bfamdef+}
\end{gather*}
for all $b = (b_i)\in k^n$.

Thus the factorization index $[L (2n+2, k): kH]^f$ is equal to the number of types of isomorphisms of Lie algebras of the set
\begin{gather*}
\{\mathfrak{l}_{(a)} (2n+1, k), \mathfrak{l}'_{(b)} (2n+1, k), \mathfrak{l}''_{(b)} (2n+1, k) \,|\, a\in k^n \setminus \{0\}, b \in k^n \}.
\end{gather*}
\end{Proposition}

\begin{proof}
$\mathfrak{l} (2n+1, k)$ is a~complement of $kH$ in $L (2n+2, k)$ and we can write $L (2n+2, k) = kH \bowtie \mathfrak{l} (2n+1,
k)$, where the bicrossed product is associated to the matched pair given in~\eqref{mpcanon}.
Hence, by~\cite[Theorem~4.3]{am-2013a} any other complement $\mathfrak{C}$ of $kH$ in $L (2n+2, k)$ is isomorphic to
an~$r$-deformation of $\mathfrak{l} (2n+1, k)$, for some deformation map $r: \mathfrak{l} (2n+1, k) \to kH$ of the matched
pair~\eqref{mpcanon}.
These are described in Lemma~\ref{defmaps}.
The Lie algebra $\mathfrak{l}_{(a)} (2n+1, k)$ is precisely the $r_a$-deformation of $\mathfrak{l} (2n+1, k)$, where $r_a$ is
given by~\eqref{def1a}.
On the other hand the $r_{(b, c)}$-deformation of $\mathfrak{l} (2n+1, k)$, where $r_{(b, c)}$ is given by~\eqref{def2a} for
some $(b = (b_i), c) \in k^n \times k$, is the Lie algebra denoted by $\mathfrak{l}_{(b, c)} (2n+1, k)$ having the bracket given
for any $i = 1, \dots, n$ by
\begin{gather*}
[E_i, F_j]_{(b, c)}:= - b_j E_i,
\qquad
[E_i, G]_{(b, c)}:= (1-c) E_i,
%\label{2afamdef+}
\\
[F_i, F_j]_{(b, c)}:= b_j F_i - b_i F_j,
\qquad
[F_i, G]_{(b, c)}:= (c-1) F_i - b_i G
%\label{2bfamdef++}
\end{gather*}
for all $(b = (b_i), c) \in k^n \times k$.
Now, for $c \neq 1$ we can see that $\mathfrak{l}_{(b, c)} (2n+1, k) \cong \mathfrak{l}'_{(b)} (2n+1, k)$ (by sending~$G$ to
$(c-1)^{-1} G$) while $\mathfrak{l}_{(b, 1)} (2n+1, k) = \mathfrak{l}''_{(b)} (2n+1, k)$ and we are done.
\end{proof}
\begin{Remark}\label{cind}
An attempt to compute $[L (2n+2, k): kH]^f$ for an arbitrary integer~$n$ is hopeless.
However, one can easily see that $\mathfrak{l}'_{(0)} (2n+1, k) = \mathfrak{l} (2n+1, k)$ and $\mathfrak{l}''_{(0)} (2n+1, k) =
k^{2n+1}_0$, the Abelian Lie algebra of dimension $2n +1$.
Thus, $[L (2n+2, k): kH]^f \geq 2$.
The case $n = 1$ is presented below.
\end{Remark}

\begin{Example}
\label{n1}
Let~$k$ be a~f\/ield of characteristic $\neq 2$ and consider $\{E, F, G\}$ the basis of $\mathfrak{l} (3, k)$ with the bracket
given by $[E, G] = E$ and $[G, F] = F$.
Then, the factorization index $[L (4, k): kH]^f = 3$.
More precisely, the isomorphism classes of all complements of $kH$ in $L (4, k)$ are represented by the following three Lie
algebras: $\mathfrak{l} (3, k)$, $k^3_0$ and the Lie algebra $L_{-1}$ having $\{E, F, G\}$ as a~basis and the bracket given~by
\begin{gather*}
[F, E] = F,
\qquad
[E, G] = - G.
\end{gather*}
Since $\charop (k) \neq 2$ the Lie algebras $\mathfrak{l} (3, k)$ and $L_{-1}$ are not isomorphic~\cite[Exercise 3.2]{EW}.
For $a\in k^*$ the Lie algebra $\mathfrak{l}_{(a)} (3, k)$ has the bracket given by $[E, F] = -a F$ and $[E, G] = -a G$.
Thus, $\mathfrak{l}_{(a)} (3, k) \cong \mathfrak{l}_{(1)} (3, k)$, and the latter is isomorphic to the Lie algebra $L_{-1}$.
On the other hand we have: $\mathfrak{l}''_{(0)} (3, k) = k^3_0$ and for $b \neq 0$ we can easily see that $\mathfrak{l}''_{(b)}
(3, k) \cong \mathfrak{l}''_{(1)} (3, k) \cong \mathfrak{l} (3, k) $.
Finally, $\mathfrak{l}'_{(0)} (3, k) = \mathfrak{l} (3, k) $ and for $b \neq 0$ we have that $\mathfrak{l}'_{(b)} (3, k) \cong
\mathfrak{l}'_{(1)} (3, k)$~-- the latter is the Lie algebra having $\{f_1, f_2, f_3 \}$ as a~basis and the bracket given by
$[f_1, f_2] = - f_1$, $[f_1, f_3] = f_1$ and $[f_3, f_2] = f_2 + f_3$.
This Lie algebra is also isomorphic to $\mathfrak{l} (3, k)$, via the isomorphism which sends $f_1$ to~$E$, $f_3$ to~$G$ and
$f_2$ to $F - G$.
\end{Example}
Let~$k$ be a~f\/ield of characteristic $\neq 2$ and $\mathfrak{l}^2 (2n+1, k)_{(A, D, \delta)}$ the Lie algebra of
Theorem~\ref{teorema11}.
In order to simplify computations we will assume $A = D:= I_n$ and $\delta:= (1, 0, \dots, 0, 1)$.
Let $\mathfrak{m} (2n+2, k):= \mathfrak{l}^2 (2n+1, k)_{(I_n, I_n, (1, 0, \dots, 0, 1))}$ be the $(2n+2)$-dimensional Lie
algebra having $\{E_i, F_i, G, H \,|\, i = 1, \dots, n \}$ as a~basis and the bracket def\/ined for any $i = 1, \dots, n$ by
\begin{gather*}
[E_i, G ] = E_i,
\qquad
[G, F_i ] = F_i,
\qquad
[E_i, H ] = E_i,
\qquad
[F_i, H ] = F_i,
\qquad
[G, H ] = E_1 + F_n.
\end{gather*}
We consider the Lie algebra extension $kH \subset \mathfrak{m} (2n+2, k)$, where $kH \cong k_0$ is the Abelian Lie algebra of
dimension~$1$.
Of course, $\mathfrak{m} (2n+2, k)$ factorizes through $kH$ and $\mathfrak{l} (2n+1, k)$, i.e.~$\mathfrak{m}(2n+2,k)=kH\bowtie\mathfrak{l} (2n+1, k)$.
Moreover, the canonical matched pair $\triangleleft: \mathfrak{l} (2n+1, k) \times kH \to \mathfrak{l} (2n+1, k)$ and
$\triangleright: \mathfrak{l} (2n+1, k) \times kH \to kH$ associated to this factorization is given as follows:
\begin{gather}
E_i \triangleleft H:= E_i,
\qquad
F_i \triangleleft H:= F_i,
\qquad
G \triangleleft H:= E_{1} + F_{n}
\label{mpcanon2}
\end{gather}
and all undef\/ined actions are zero.
In particular, we should notice that the left action $\triangleright: \mathfrak{l} (2n+1, k) \times kH \to kH$ is trivial.
Next, we describe the set ${\mathcal D}{\mathcal M}(\mathfrak{l} (2n+1, k), kH | (\triangleright, \triangleleft))$
of all deformation maps of the matched pair $(kH, \mathfrak{l} (2n+1, k), \triangleright, \triangleleft)$ given
by~\eqref{mpcanon2}.

\begin{Lemma}
\label{defmaps2}
Let~$k$ be a~field of characteristic $\neq 2$.
Then there exists a~bijection
\begin{gather*}
{\mathcal D}{\mathcal M} \big(\mathfrak{l} (2n+1, k), kH \big| (\triangleright, \triangleleft) \big) \cong \big(k^n \setminus
\{0\} \big) \sqcup \big(k^n \setminus \{0\} \big) \sqcup k.
\end{gather*}
The bijection is given such that the deformation map $r = r_a: \mathfrak{l} (2n+1, k) \to kH$ associated to $a = (a_i) \in k^n
\setminus \{0\}$ is given~by
\begin{gather}
\label{def1a1}
r (E_i):= a_i H,
\qquad
r(F_i):= 0,
\qquad
r (G):= (a_{1} - 1)H
\end{gather}
the deformation map $r = r_b: \mathfrak{l} (2n+1, k) \to kH$ associated to another $b = (b_i) \in k^n \setminus \{0\}$ is given
by
\begin{gather}
\label{def1a2}
r (E_i):= 0,
\qquad
r(F_i):= b_{i} H,
\qquad
r (G):= (b_{n} + 1)H,
\end{gather}
while the deformation map $r = r_{c}: \mathfrak{l} (2n+1, k) \to kH$ associated to $c \in k$ is given~by
\begin{gather}
\label{def2a3}
r (E_i):= 0,
\qquad
r(F_i):= 0,
\qquad
r (G):= c H
\end{gather}
for all $i = 1, \dots, n$.
\end{Lemma}

\begin{proof}
Any linear map $r: \mathfrak{l} (2n+1, k) \to kH$ is uniquely determined by a~triple $(a = (a_i)$, $b = (b_i), c) \in k^n \times
k^n \times k$ via: $r(E_i):= a_i H$, $r(F_i):= b_i H$ and $r(G):= c H $, for all $i = 1, \dots, n$.
We only need to check when such a~map $r = r_{(a, b, c)}$ is a~deformation map.
Since $kH$ is the Abelian Lie algebra and the left action $\triangleright: \mathfrak{l} (2n+1, k) \times kH \to kH$ is trivial,
equation~\eqref{factLie} comes down~to
\begin{gather}
\label{defabc2}
r([x, y]) = r \big(y \triangleleft r(x) - x \triangleleft r(y) \big).
\end{gather}
Since~\eqref{defabc2} is symmetrical it is enough to check it only for pairs of the form $(E_{i}, E_{j})$, $(F_{i}, F_{j})$,
$(E_{i}, F_{j})$, $(E_{i}, G)$, and $(F_{i}, G)$, for all~$i, j = 1, \dots, n$.
It is straightforward to see that~\eqref{defabc2} is trivially fulf\/illed for the pairs $(E_{i}, E_{j})$, $(F_{i}, F_{j})$ and
$(E_{i}, F_{j})$.
Moreover,~\eqref{defabc2} evaluated for $(E_{i}, G)$ and respectively $(F_{i}, G)$ yields $a_{i}(a_{1} + b_{n} - c - 1) = 0$ and
$b_{i}(a_{1} + b_{n} - c + 1) = 0$ for all $i = 1, \dots, n$.
Therefore, keeping in mind that we work over a~f\/ield of characteristic $\neq 2$, the triples $(a = (a_i), b = (b_i), c) \in k^n
\times k^n \times k$ for which $r_{(a, b, c)}$ becomes a~deformation map are given as follows: $(a = (a_i) \in k^{n}\setminus
\{0\}, b = 0, c = a_{1} - 1)$, $(a = 0, b = (b_i) \in k^{n}\setminus \{0\}, c = b_{n} + 1)$ and $(a = 0, b = 0, c \in k)$.
The corresponding deformation maps are exactly those listed above.
\end{proof}

The next result describes all deformations of $\mathfrak{l} (2n+1, k)$ associated to the canonical matched pair $(kH,
\mathfrak{l} (2n+1, k), \triangleright, \triangleleft)$ given by~\eqref{mpcanon2}.

\begin{Proposition}
\label{mpdef2}
Let~$k$ be a~field of characteristic $\neq 2$ and the extension of Lie algebras $kH \subset \mathfrak{m} (2n+2, k)$.
Then a~Lie algebra $\mathfrak{C}$ is a~complement of $kH$ in $\mathfrak{m} (2n+2, k)$ if and only if $\mathfrak{C}$ is
isomorphic to one of the Lie algebras from the three families defined below:

$\bullet$ the Lie algebra $\overline{\mathfrak{l}}_{(a)} (2n+1, k)$ having the bracket def\/ined for any $i = 1, \dots, n$ by
\begin{gather*}
[E_{i}, E_{j}]_{a}:= a_{j} E_{i} - a_{i} E_{j},
\qquad
[E_i, F_j]_{a}:= - a_i F_j,
\\
[E_i, G]_{a}:= a_{1} E_{i} - a_{i} (E_{1} + F_{n}),
\qquad
[G, F_i]_{a}:= (2-a_{1}) F_i
%\label{Lie111}
\end{gather*}
for all $a = (a_i) \in k^n \setminus \{0\}$.

$\bullet$ the Lie algebra $\overline{\mathfrak{l'}}_{(b)} (2n+1, k)$ having the bracket def\/ined for any $i = 1, \dots, n$ by
\begin{gather*}
[F_{i}, F_{j}]_{b}:= b_{j} F_{i} - b_{i} F_{j},
\qquad
[E_i, F_j]_{b}:= b_j E_i,
\\
[E_i, G]_{b}:= (2 + b_n) E_i,
\qquad
[G, F_i]_{b}:= b_{i}(E_{1} + F_{n}) - b_{n} F_{i}
%\label{Lie112}
\end{gather*}
for all $b = (b_i)\in k^n \setminus \{0\}$.

$\bullet$ the Lie algebra $\overline{\mathfrak{l''}}_{(c)} (2n+1, k)$ having the bracket def\/ined for any $i = 1, \dots, n$ by
\begin{gather*}
%\label{Lie113}
[E_i, G]_c:= (1+c) E_i,
\qquad
[G, F_i]_c:= (1-c) F_i
\end{gather*}
for all $c \in k$.

Thus the factorization index $[\mathfrak{m} (2n+2, k): kH]^f$ is equal to the number of types of isomorphisms of Lie algebras of
the set
\begin{gather*}
\big\{\overline{\mathfrak{l}}_{(a)} (2n+1, k), \,\overline{\mathfrak{l}'}_{(b)} (2n+1, k), \, \overline{\mathfrak{l}''}_{(c)} (2n+1, k) \,|\,
a, b \in k^n \setminus \{0\},\, c \in k \big\}.
\end{gather*}
\end{Proposition}

\begin{proof}
As in the proof of Proposition~\ref{mpdef} we make use of~\cite[Theorem~4.3]{am-2013a}.
More precisely, this implies that all complements $\mathfrak{C}$ of $kH$ in $\mathfrak{m} (2n+2, k)$ are isomorphic to
an~$r$-deformation of $\mathfrak{l} (2n+1, k)$, for some deformation map $r: \mathfrak{l} (2n+1, k) \to kH$ of the matched
pair~\eqref{mpcanon2}.
These are described in Lemma~\ref{defmaps2}.
By a~straightforward computation it can be seen that $\overline{\mathfrak{l}}_{(a)} (2n+1, k)$ is exactly the complement
corresponding to the deformation map given by~\eqref{def1a1}, $\overline{\mathfrak{l'}}_{(b)} (2n+1, k)$ corresponds to the
deformation map given by~\eqref{def1a2} while $\overline{\mathfrak{l''}}_{(c)} (2n+1, k)$ is implemented by the deformation map
given by~\eqref{def2a3}.
\end{proof}

\begin{Example}
%\label{n2}
Let~$k$ be a~f\/ield of characteristic $\neq 2$.
Then, the factorization index $[\mathfrak{m} (4, k): kH]^f$ depends essentially on the f\/ield~$k$.
We will prove that all complements of $kH$ in $\mathfrak{m} (4, k)$ are isomorphic to a~Lie algebra of the form:
\begin{gather*}
L_{\alpha}: \ [x, z] = x,
\qquad
[y, z] = \alpha y,
\qquad
\text{with}
\quad
\alpha \in k.
\end{gather*}
Hence, $[\mathfrak{m} (4, k): kH]^f = \infty $, if $|k| = \infty $ and $[\mathfrak{m} (4, k): kH]^f = (1 + p^n)/2$, if $|k| =
p^n$, where $p\geq 3$ is a~prime number.
Indeed, for $n=1$, the Lie algebras described in Proposition~\ref{mpdef2} become
\begin{gather*}
\overline{\mathfrak{l}}_{(a)} (3, k): \quad  [E, F]_{a}:= - aF, \qquad [E, G]_{a}:= - aF, \qquad [G, F]_{a}:= (2-a) F,
\\
\overline{\mathfrak{l'}}_{(b)} (3, k): \quad  [E, F]_{b}:= b E, \qquad [E, G]_{b}:= (2 + b) E, \qquad [G, F]_{b}:= b E,
\\
\overline{\mathfrak{l''}}_{(c)} (3, k): \quad  [E, G]_c:= (1+c) E, \qquad [G, F]_c:= (1-c) F,
\end{gather*}
$a, b \in k^{*}$, $c \in k$.
To start with, we should notice that the f\/irst two Lie algebras $\overline{\mathfrak{l}}_{(a)}$ and
$\overline{\mathfrak{l'}}_{(b)}$ are isomorphic for all~$a, b \in k^{*}$.
The isomorphism $\gamma: \overline{\mathfrak{l}}_{(a)} \to \overline{\mathfrak{l'}}_{(b)}$ is given as follows
\begin{gather*}
\begin{split}
& \gamma(E)  :=  2^{-1}(b-a) E + 2^{-1} (b-a+2) F + 2^{-1}(a-b) G,
\qquad
\gamma(F):= E,
\\
& \gamma(G)  :=  2^{-1} (b-a+4) E + 2^{-1}(b-a+4) F + 2^{-1}(a-b-2) G.
\end{split}
\end{gather*}
Moreover, the map $\varphi: \overline{\mathfrak{l}}_{(a)} \to L_{0}$ given by
\begin{gather*}
\varphi(E):= y + az,
\qquad
\varphi(F):= x,
\qquad
\varphi(G):= x + y + (a-2) z
\end{gather*}
is an isomorphism of Lie algebras for all $a \in k^{*}$.
Therefore, the f\/irst two Lie algebras are both isomorphic to $L_{0}$ for all~$a, b \in k^{*}$.
We are left to study the family $\overline{\mathfrak{l''}}_{(c)}$.
If $c = -1$ then $\overline{\mathfrak{l''}}_{(-1)}$ is again isomorphic to $L_{0}$.
Suppose now that $c \neq -1$.
Then the map $\psi: \overline{\mathfrak{l''}}_{(c)} \to L_{(c-1)(c+1)^{-1}}$ given by
\begin{gather*}
\psi(E):= x,
\qquad
\psi(F):= y,
\qquad
\psi(G):= (c+1) z
\end{gather*}
is an isomorphism of Lie algebras.
Finally, we point out here that if $\alpha \notin \{\beta, \beta^{-1} \}$ then $L_{\alpha}$ is not isomorphic to $L_{\beta}$
(see, for instance~\cite[Exercise 3.2]{EW}) and the conclusion follows.
\end{Example}

\begin{Remark}\label{2noi}
We end this section with two more applications.
The deformation of a~given Lie algebra $\mathfrak{h}$ associated to a~matched pair $(\mathfrak{g}, \mathfrak{h}, \triangleright,
\triangleleft)$ of Lie algebras and to a~deformation map~$r$ as def\/ined by~\eqref{rLiedef} is a~very general method of
constructing new Lie algebras out of a~given Lie algebra.
It is therefore natural to ask if the properties of a~Lie algebra are preserved by this new type of deformation.
We will see that in general the answer is negative.
First of all we remark that the Lie algebra $\mathfrak{h}:= \mathfrak{l} (2n+1, k)$ is metabelian,
that is $[[\mathfrak{h},\mathfrak{h}], [\mathfrak{h},\mathfrak{h}]] = 0$.
Now, if we look at the matched pair deformation $\mathfrak{h}_r = \mathfrak{l}_{(a)} (2n+1, k)$ of $\mathfrak{h}$ given
by~\eqref{1famdef} of Proposition~\ref{mpdef}, for $a = (a_i) \in k^n \setminus \{0\}$ we can easily see that
$\mathfrak{l}_{(a)} (2n+1, k)$ is not a~metabelian Lie algebra, but a~$3$-step solvable Lie algebra.
Thus the property of being metabelian is not preserved by the~$r$-deformation of a~Lie algebra.

Next we consider an example of a~somewhat dif\/ferent nature.
First recall~\cite{medina} that a~Lie algebra $\mathfrak{h}$ is called \emph{self-dual} (or \emph{metric}) if there exists
a~non-degenerate invariant bilinear form $B: \mathfrak{h} \times \mathfrak{h} \to k$, i.e.~$B([a, b], c) = B(a, [b, c])$, for
all~$a,b,c\in \mathfrak{h}$.
Self-dual Lie algebras generalize f\/inite-dimensional complex semisimple Lie algebras (the second Cartan's criterion shows that
any f\/inite-dimensional complex semisimple Lie algebra is self-dual since its Killing form is non-degenerate and invariant).
Besides the mathematical interest in studying self-dual Lie algebras, they are also important and have been intensively studied
in physics~\cite{fig, pelc}.
Now, $\mathfrak{h}:= \mathfrak{l} (2n+1, k)$ is not a~self-dual Lie algebra since if $B: \mathfrak{l} (2n+1, k) \times
\mathfrak{l} (2n+1, k) \to k$ is an arbitrary invariant bilinear form then we can easily prove that $B(E_i, -) = 0$ and thus any
invariant form is degenerate.
On the other hand, the~$r$-deformation of $\mathfrak{l} (2n+1, k)$ denoted by $\mathfrak{l}_{(0)}'' (2n+1, k)$ in
Remark~\ref{cind} is self-dual since it is just the $(2n+1)$-dimensional Abelian Lie algebra.
\end{Remark}

\section{Two open questions}

The paper is devoted to the factorization problem and its converse, the classifying complements problem, at the level of Lie
algebras.
Both problems are very dif\/f\/icult ones; even the case considered in this paper, namely $\mathfrak{g} = k_0$, illustrates the
complexity of the two problems.
We end the paper with the following two open questions:

\textbf{Question 1.} \emph{Let $n \geq 2$.
Does there exist a~Lie algebra $\mathfrak{h}$ and a~matched pair of Lie algebras $(\mathfrak{gl} (n, k), \mathfrak{h},
\triangleleft, \triangleright)$ such that $\mathfrak{gl} (n, k) \bowtie \mathfrak{h} \cong \mathfrak{gl} (n+1, k)$?}

A more restricted version of this question is the following: does the canonical inclusion $\mathfrak{gl} (n, k) \hookrightarrow
\mathfrak{gl} (n +1, k)$ have a~complement that is a~Lie subalgebra of $\mathfrak{gl} (n +1, k)$? Although it seems unlikely for
such a~complement to exist we could not f\/ind any proof or reference to this problem in the literature.

Secondly, having~\cite[Corollary 3.2]{am-2012} as a~source of inspiration we ask:

\textbf{Question 2.} \emph{Let $n \geq 2$.
Does there exist a~matched pair of Lie algebras $(\mathfrak{g}, \mathfrak{h}, \triangleleft, \triangleright)$ such that
any~$n$-dimensional Lie algebra $\mathfrak{L}$ is isomorphic to an~$r$-deformation of $\mathfrak{h}$ associated to this matched
pair?}

At the level of groups, question 2 has a~positive answer by considering the canonical matched pair associated to the
factorization of $S_{n+1}$ by $S_n$ and the cyclic group $C_n$.

\subsection*{Acknowledgements}

We would like to thank the referees for their comments and suggestions that substantially improved the f\/irst version of this
paper.
A.L.~Agore is research fellow `Aspirant' of FWO-Vlaanderen.
This work was supported by a~grant of the Romanian National Authority for Scientif\/ic Research, CNCS-UEFISCDI, grant
no.~88/05.10.2011.

\pdfbookmark[1]{References}{ref}
\LastPageEnding

\end{document}